\documentclass[12pt]{article}

\usepackage{amsmath,amsfonts,upref}
\usepackage{mylabels}
\usepackage{mymargs}
\usepackage{graphics}

\newcommand{\nz}{{n_z}}
\newcommand{\hpipi}{(\pi/2,\pi)}
\newcommand{\vbf}{\mathbf{v}}
\newcommand{\zbf}{\mathbf{z}}

\newcommand{\Zbb}{\mathbb{Z}}
\newcommand{\Rbb}{\mathbb{R}}
\newcommand{\Nbb}{\mathbb{N}}
\newcommand{\abs}[1]{\lvert#1\rvert}

\newtheorem{theorem}{Theorem}[section]

\begin{document}
\sloppy

\pagestyle{myheadings}

\title{\bf Orthonormal Compactly Supported Wavelets with Optimal Sobolev Regularity}

\author{Harri Ojanen%
\thanks{The author was supported in part by the Academy of Finland and
    the Finnish Society of Sciences and Letters.}}
    
\date{June 30, 1998}

\maketitle

{\footnotesize
\begin{quote}
{\bf Abstract}: Numerical optimization is used to construct new
orthonormal compactly supported wavelets with Sobolev regularity
exponent as high as possible among those mother wavelets with a fixed
support length and a fixed number of vanishing moments. The increased
regularity is obtained by optimizing the locations of the roots the
scaling filter has on the interval $\hpipi$. The results
improve those obtained by I.~Daubechies [Comm.\ Pure Appl.\ Math.\
\textbf{41} (1988), 909--996], H.~Volkmer [SIAM J.\ Math.\ Anal.\
\textbf{26} (1995), 1075--1087], and P.~G.~Lemari\'e-Rieusset and
E.~Zahrouni [Appl.\ Comput.\ Harmon.\ Anal.\ \textbf{5} (1998),
92--105].

\medbreak



{\bf E-mail address}: ojanen@math.rutgers.edu

{\bf Address}: Department of Mathematics - Hill Center;
    Rutgers, the State University of New Jersey;
    110~Frelinghuysen Rd; 	
    Piscataway NJ 08854-8019; USA
\end{quote}
}

\markright{Wavelets with optimal Sobolev regularity}

\thispagestyle{empty}

\section{Introduction}\mylabel{s:intro}

To construct orthonormal compactly supported wavelets we use as the
starting point the dilation equation
\begin{equation}\mylabel{e:dileqn}
    \phi(x) = \sqrt2 \sum_{k\in\Zbb} 
	c_k \phi(2x-k), \qquad x\in\Rbb, \qquad c_k\in\Rbb,
\end{equation}
for the scaling function or father function~$\phi$. A solution exists
in the sense of distributions when $\sum_k c_k=\sqrt2$ and it has
compact support if only finitely many of the filter coefficients $c_k$
are non-zero (in fact, if $c_k\not=0$ only for $k\in\{0,\dots,n\}$
then the support of~$\phi$ is contained in~$[0,n]$). The solution of
the dilation equation is given by
\begin{equation}\mylabel{e:m0prod}
    \widehat\phi(\xi) = \prod_{j=1}^\infty m_0(2^{-j}\xi),
\end{equation}
where the trigonometric polynomial $m_0(\xi)$ is the scaling filter
of~$\phi$,
\begin{equation}\mylabel{e:m0}
    m_0(\xi) = \frac1{\sqrt2} \sum_{k\in\Zbb} c_k e^{-ik\xi}.
\end{equation}

The condition when the solution to~\eqref{e:dileqn} yields wavelets can
be expressed in terms of the scaling filter $m_0(\xi)$.
If $m_0$ satisfies
\begin{equation}\mylabel{e:ortho}
    \abs{m_0(\xi)}^2 + \abs{m_0(\pi+\xi)}^2 = 1, \qquad \xi\in\Rbb,
\end{equation}
$m_0(0)=1$, and a technical condition
called the Cohen criterion holds, then there exists a scaling function and
wavelet pair \cite{Cohe90}, \cite[Section~6.3]{Daub92}. We will not
need exact details about Cohen's condition, we only 
use that if $m_0$ does not vanish on
$[-\pi/2,\pi/2]$, then it satisfies the criterion.

If $\phi$ is such that it gives rise to
a wavelet, the wavelet $\psi$ is obtained by
\[
    \psi(x) = \sum_{k\in\Zbb}  (-1)^k c_{1-k} \phi(2x-k).
\]
Note that then $\phi$ and $\psi$ have the same regularity
properties.

The original construction of orthonormal compactly supported wavelets
by I.~Daubechies~\cite{Daub88} considers a trigonometric polynomial
$m_0$ with at most $2N$ non-zero coefficients $c_0$, \dots,
$c_{2N-1}$. The polynomial $m_0$ is required to
satisfy~\eqref{e:ortho} and to have a zero at $\pi$ of maximal
order~$N$ (the order of the zero is also the number of vanishing
moments the corresponding wavelet has). The resulting equations can be
solved and the constructed $m_0$ has no roots on $(-\pi,\pi)$, hence it
satisfies the Cohen criterion.  (\cite{Daub92} is a comprehensive
reference for these results.)

To construct smoother wavelets our approach is to reduce the order of
the zero $m_0$ has at $\pi$ and instead to introduce zeros on the
interval $\hpipi$. Keeping $N$ and the number of zeros on $\hpipi$ fixed
we use numerical optimization to choose the locations of the roots so
that the resulting wavelets are as smooth as possible. We are able to
construct wavelets that are more regular than those introduced in
\cite{Daub88},~\cite{Lema98}, and~\cite{Volk95}.

This paper is organized as follows: in Section~\ref{s:res} we state
our results and compare them to the above mentioned papers. The
methods used are described in Section~\ref{s:methods}, and in
Section~\ref{s:optim} the optimization is discussed in detail. We
comment in Section~\ref{s:other-param} on some other approaches that
were tried. Finally, filter coefficients for selected most regular
wavelets are listed in section~\ref{s:coef}.

\section{Results}\mylabel{s:res}

In the numerical optimization for smoothness we studied wavelets with
up to $40$ filter coefficients and whose scaling filter has up to four
roots on the interval $\hpipi$. The Sobolev regularity exponents of
the most regular wavelets found are shown in Table~\ref{t:s0res}. (By
the Sobolev regularity exponent $s_0$ of $\phi$ we mean $s_0 =
\sup\{s\colon \phi\in H^s\}$, where $H^s$ is the Sobolev space
consisting of all $f\in L^2(\Rbb)$ such that $(1+\abs{\xi}^2)^{s/2}
\widehat{f}(\xi)\in L^2(\Rbb)$.)  The results are also summarized in
Figure~\ref{f:s0res}.

\begin{table}[t]
\begin{center}
\begin{tabular}{ccccccc}
&& \multicolumn{5}{c}{Sobolev regularity exponent $s_0$} \\
$2N$ &\quad& $\nz=0$ & $\nz=1$ & $\nz=2$ & $\nz=3$ & $\nz=4$ \\[.5ex]\hline
\\[-2ex]
    2   &&      0.50  \\
    4   &&      1.00  \\
    6   &&      1.42  &    1.00   \\
    8   &&      1.78  &    1.82   \\
    10  &&      2.10  &    2.26   &    1.00  \\
    12  &&      2.39  &    2.66   &    2.00  \\
    14  &&      2.66  &    3.02   &    3.00  &    1.00   \\
    16  &&      2.91  &    3.37   &    3.48  &    2.00   \\
    18  &&      3.16  &    3.72   &    3.92  &    3.00   &    1.00   \\
    20  &&      3.40  &    4.07   &    4.32  &    4.00   &    2.00   \\
    22  &&      3.64  &    4.42   &    4.73  &    4.76   &    3.00   \\
    24  &&      3.87  &    4.78   &    5.14  &    5.22   &    4.00   \\
    26  &&      4.11  &    5.14   &    5.55  &    5.67   &    5.00   \\
    28  &&      4.34  &    5.50   &    5.95  &    6.10   &    5.84   \\
    30  &&      4.57  &    5.85   &    6.34  &    6.51   &    6.38   \\
    32  &&      4.79  &    6.19   &    6.71  &    6.89   &    6.88   \\
    34  &&      5.02  &    6.52   &    7.09  &    7.25   &    7.30   \\
    36  &&      5.24  &    6.83   &    7.45  &    7.62   &    7.69   \\
    38  &&      5.47  &    7.15   &    7.81  &    7.99   &    8.08   \\
    40  &&      5.69  &    7.46   &    8.17  &    8.37   &    8.51   
\end{tabular}
\vspace*{.5ex}
\end{center}
\caption{The best Sobolev regularity exponents found. The length of
the filter is $2N$ and the filter $m_0$ has $n_z$ roots on the
interval $\hpipi$. For comparison the column $\nz=0$ gives the Sobolev
exponents of the Daubechies wavelets~\cite{Daub88}, the other columns 
are new.}\mylabel{t:s0res}
\end{table}

The column $\nz=0$ gives the Sobolev exponents for the Daubechies
wavelets constructed in~\cite{Daub88}. For $N\geq4$ we have 
constructed wavelets which are more regular. For large values of $N$
there are four families of wavelets corresponding to different numbers of
roots on $\hpipi$, and the regularity increases with the number of roots.

To relate these results to H\"older regularity, it can be shown that the
H\"older regularity exponent $\alpha_0=\sup \{\alpha\colon
\phi\in C^\alpha\}$ satisfies $\alpha_0 \in [s_0-1/2,s_0]$
\cite[Corollary~9.9]{Vill92}.  Here $C^\alpha$ consists of functions
that are continuously differentiable up to order equal to the integer
part of $\alpha$, and the highest order continuous derivative 
belongs to the (standard) H\"older class with exponent equal to the
fractional part of~$\alpha$.

Note that in each column $\nz=1$, \dots $4$, the first few entries are
smaller than in the previous columns on the same row. This is not
surprising, since it can be shown that in order to have $\psi\in H^n$,
the wavelet $\psi$ must have at least $n+1$ vanishing
moments~\cite[Proposition~8.3]{Vill92}. The wavelets corresponding to
the entries in the table have $N-2\nz$ vanishing moments.

We show below (see Section~\ref{s:dof2}) that the entries in the column
${\nz=1}$ are best possible in the sense that the corresponding wavelets are
smoothest possible among those with $2N\geq8$ coefficients and at least
$N-2$ vanishing moments.

\begin{figure}
\begin{center}
\scalebox{0.6}{\includegraphics{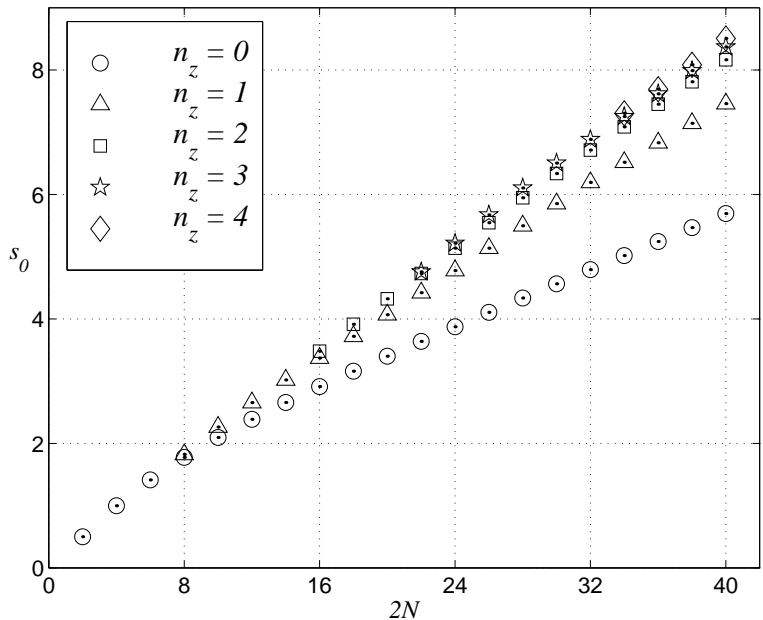}}
\end{center}
\caption{The Sobolev regularity exponents $s_0$ from Table~\ref{t:s0res}.
For clarity only
those are shown that are better than the ones with the same length
filter but with fewer roots on $\hpipi$.}\mylabel{f:s0res}
\end{figure}

For higher values of $\nz$ the exponents listed in Table~\ref{t:s0res}
are best possible among those wavelets with $2N$ coefficients and
whose scaling filter $m_0$ has the indicated number of zeros on
$\hpipi$. Some of these were subjected to extensive tests by choosing
many different initial approximations for the optimizer. Note also how
the results, except for the first few entries, lie on almost straight
lines in Figure~\ref{f:s0res}. See also Figures~\ref{f:s0z1}
and~\ref{f:s0z2} which show typical examples of the regularity
exponent as a function of the location of the roots.

However, it is not completely clear if the results for $n_z\geq2$ are
best possible if we consider wavelets with the same number of
vanishing moments as above but do not insist on $m_0$ having roots on
$\hpipi$. We come back to this at the end of Section~\ref{s:dofn}.

We remark that the case $2N=10$, $\nz=1$ appears already in~\cite{Ojan91}.

More regular wavelets than the original Daubechies family
have been constructed by I.~Daubechies~\cite{Daub93},
H.~Volkmer~\cite{Volk95},
P.~G.~Lemari\'e-Rieusset and E.~Zahrouni~\cite{Lema98}.
We discuss the wavelets from~\cite{Daub93} in Section~\ref{s:dof1}.

H.~Volkmer~\cite{Volk95} computes only the asymptotic regularity ratio
$\lim_{N\rightarrow\infty}\sigma/2N$ and constructs wavelets for which
the ratio is close to $0.19$. Note
that for many entries in Table~\ref{t:s0res} the ratio
$\sigma/2N$ is over $0.21$. 

P.~G.~Lemari\'e-Rieusset and E.~Zahrouni~\cite{Lema98} discuss what they
call the Matzinger wavelets: the scaling filter has a multiple root at
$2\pi/3$ and no other roots on $\hpipi$. 
Our results are better, e.g., the smoothest
one they construct has $34$ coefficients and Sobolev exponent $6.13$.
Table~\ref{t:s0res} shows that already with $\nz=1$ our wavelet with the
same number of coefficients is more regular. In fact, for each value of
$N$ the $\nz=1$ wavelet listed in Table~\ref{t:s0res} is better than the
smoothest wavelet discussed in~\cite{Lema98}.

\section{Methods}\mylabel{s:methods}

We begin by describing the theorem that allows us to compute exactly the
Sobolev regularity exponent of solutions to~\eqref{e:dileqn}.
Define another trigonometric polynomial, $r(\xi)$, by factoring out all
zeros of $\abs{m_0(\xi)}^2$ at $\xi=\pi$, that is, let
\begin{equation}\mylabel{e:m0def}
    \abs{m_0(\xi)}^2 = \left( \frac{1+\cos\xi}{2} \right)^M r(\xi),
\end{equation}
where $M\in\Nbb$ is chosen so that $r(\pi)\not=0$. The function $r$ is
used to define an operator on $C[0,\pi]$ by
\begin{equation}\mylabel{e:Trdef}
    T_r u(\xi) = r(\xi/2) u(\xi/2) + r(\pi-\xi/2) u(\pi-\xi/2).
\end{equation}
The Sobolev regularity exponent of a solution to the dilation
equation~\eqref{e:dileqn} is given by the following theorem, proved
independently by L.~Villemoes and T.~Eirola:

\begin{theorem}[\cite{Eiro92,Vill92}]\mylabel{t:sobexp} Suppose $\phi$ is
a solution to~\eqref{e:dileqn}, $M$, $r$, and $T_r$ are as above, and
$r$ satisfies the Cohen criterion. Then $\phi\in H^s$ if and only if
$s<s_0=M-\log_4 \rho(T_r)$, where $\rho(T_r)$ is the spectral radius of
$T_r$ on $C[0,\pi]$.
\end{theorem}

If $r$ is a polynomial of cosines of degree at most $d$, i.e.,
$r(\xi)=\sum_{k=0}^d b_k \cos(k\xi)$, then $T_r$ maps the space spanned
by $\{\cos(k\xi)\}_{k=0}^d$ into itself and the spectral radius of $T_r$ can be
calculated from the eigenvalues of a finite matrix~\cite{Eiro92,Vill92}.

Suppose we want to look at a filter of length $2N$ (i.e., a filter
$m_0$ with at most $2N$ non-zero coefficients), such that $m_0$ has a
zero of order $M$ at $\pi$ and $\nz$ zeros on the interval $\hpipi$.
Note that $M$ is also the number of vanishing moments the associated
wavelet has (if it exists).

Define a new sequence of coefficients by
\begin{equation}\mylabel{e:ak}
    \abs{m_0(\xi)}^2 = \sum_{k=0}^{2N-1} a_k \cos(k\xi).
\end{equation}
The orthonormality condition \eqref{e:ortho} is equivalent to
\[
    a_0 = \frac12, \qquad a_{2k}=0, \quad
        \text{for $k=1$, \dots, $N-1$.}
\]
We first use the zero of $m_0$ at $\pi$ to find (linear) relations among the
$\{a_k\}$: we express
$a_k$, $k=1$, $3$, \dots, $2M-1$, in terms of $a_k$,
$k=2M+1$, $2M+3$, \dots, $2N-1$. Note that there are $N-M=2n_z$
independent parameters.  
For convenience we switch to a new set of parameters
$\vbf=(v_1,\dots,v_{n_z})$ that are simply obtained from
$a_k$, $k=2M+1$, $2M+3$, \dots, $2N-1$ by scaling the latter so that
$\vbf=(1,\dots,1)$ corresponds to the Daubechies wavelets (i.e., all
roots of $m_0$ are at $\pi$).
The function $r$ is then obtained by
factoring $\big((1+\cos\xi)/2\big)^M$ from $\abs{m_0(\xi)}^2$.
(The computations are done quickly with a symbolic algebra
system such as Mathematica by using the Taylor polynomial of~$m_0$ at~$\pi$.)

When $r$ has low degree it is easy to compute the matrix of the
operator $T_r$ on the span of $\{\cos(k\xi)\}_{k=0}^{\deg r}$
algebraically---an example is done in Section~\ref{s:dof1}. 
When the degree is high it is more convenient to
compute the matrix numerically by using the fast Fourier transform.
Finally, $s_0$ is computed from the largest
eigenvalue of this matrix.

Note that fixing the value of $\vbf$ does not uniquely determine the
wavelet, since it only determines the absolute value of~$m_0(\xi)$. For
our purposes this is enough, since all those wavelets have the same
regularity by Theorem~\ref{t:sobexp}.

\section{Optimization}\mylabel{s:optim}

\subsection{One degree of freedom}\mylabel{s:dof1}

First consider wavelets with four filter coefficients, i.e., $N=2$, and one
vanishing moment. Then
\[
    \abs{m_0(\xi)}^2 = \frac12 + a_1 \cos\xi + a_3 \cos3\xi
\]
and $m(\pi)=0$ gives $a_1=1/2-a_3$.  To compute $r$ we get
\[
    \abs{m_0(\xi)}^2 = 
	\frac12
        (1+\cos\xi) (8a_3 \cos^2\xi-8a_3 \cos\xi + 1) =
	\frac12
        (1+\cos\xi) r(\xi).    
\]
The Daubechies wavelet has two vanishing moments, hence $r(\pi)=0$,
which implies $a_3=-1/16$. Introduce the new parameter
$v$ by $a_3=-v/16$. Then
\[
    r(\xi) = -\frac v2 \cos^2\xi + \frac v2 \cos\xi + 1,
\]
and the matrix of the operator $T_r$ on the span of
$\{1,\cos\xi,\cos2\xi\}$ is easily computed to be
\[
    \frac1{4}
    \begin{pmatrix}
        8-2v    &   -2v     &   0       \\
        2v      &   2v      &   0       \\
        -v      &   8-2v    &   -v  
    \end{pmatrix}.
\]
The eigenvalues of this matrix are
\[
    \lambda_1 = 1+\sqrt{1-v}, \quad
    \lambda_2 = 1-\sqrt{1-v}, \quad \text{and} \quad
    \lambda_3 = -\frac{v}{4}.
\]
Recall that $r(\xi)$ has to be non-negative. In particular $r(\pi)\geq0$
gives $v\leq1$. Then $\lambda_1\geq1$ and $\rho(T_r)\geq1$ and so
$s_0\leq1$ for all admissible values of $v$. The Daubechies wavelet
corresponds to $v=1$ and $s_0=1$, for all other wavelets in this family
we have $v<1$, and therefore $s_0<1$.

Numerical results show this is also the case for wavelets with larger
values of $N$ and with one degree of freedom. Thus the original
Daubechies~\cite{Daub88} wavelets are the smoothest possible among
wavelets with $2N$ coefficients and at least $N-1$ vanishing moments.

I.~Daubechies~\cite{Daub93} has also constructed wavelets that have
higher H\"older regularity exponents than the ones originally
constructed in~\cite{Daub88}. In~\cite{Daub93} the cases $2N=4$ and
$2N=6$ were studied. Table~\ref{t:hold-sob1} gives both H\"older and
Sobolev exponents for two wavelets from the original construction
in~\cite{Daub88}. (The Sobolev exponents are
as in Table~\ref{t:s0res}, the H\"older exponents are
from~\cite{Daub92a}.)
Table~\ref{t:hold-sob2} gives the
exponents for the two wavelets constructed in~\cite{Daub93}. The
Sobolev exponents for the new wavelets are easily computed from
Theorem~\ref{t:sobexp}, e.g., for $2N=4$ we get
$s_0=\log_4(10/3)\approx0.87$.

Comparing the numbers in Tables~\ref{t:hold-sob1}
and~\ref{t:hold-sob2} an interesting phenomenon appears: both for
$2N=4$ and for $2N=6$ we have an example of two functions, one of which is
slightly more regular when regularity is measured in the H\"older
sense, but the opposite is true if Sobolev regularity is used.

\begin{table}
\begin{center}
\begin{tabular}{cp{2mm}cc}
    $2N$ && H\"older & Sobolev \\[.5ex]\hline \\[-2ex]
        4 && 0.55 & 1.00 \\
        6 && 1.09 & 1.42
\end{tabular}
\end{center}
\caption{H\"older and Sobolev regularity exponents for the
wavelets constructed in~\cite{Daub88}.}\mylabel{t:hold-sob1}
\end{table}

\begin{table}
\begin{center}
\begin{tabular}{cp{2mm}cc}
    $2N$ && H\"older & Sobolev \\[.5ex]\hline \\[-2ex]
        4 && 0.59 & 0.87 \\
        6 && 1.40 & 1.31
\end{tabular}
\end{center}
\caption{H\"older and Sobolev regularity exponents for the H\"older
smoothest wavelets constructed in~\cite{Daub93}. (The H\"older
exponents are lower bounds for the actual value.)}\mylabel{t:hold-sob2}
\end{table}

\subsection{Two degrees of freedom}\mylabel{s:dof2}

Optimization was done on the $\vbf=(v_1,v_2)$ parameters. Since some
values of $\vbf$ yield an $r_\vbf$ that takes on also negative values,
this is now a constrained optimization problem with the nonlinear
constraint $\min_{0\leq\xi\leq\pi} r_\vbf(\xi) \geq0$.

Figure~\ref{f:s0102} shows a plot of the Sobolev exponent as a
function of the parameters $v_1$,~$v_2$. A~typical example of the
exponent on the boundary of the feasible region is given in
Figure~\ref{f:s0z1}, where the exponent is plotted as a function of
the double root that $r$ has on $\hpipi$. Figure~\ref{f:r322} shows
a plot of $r(\xi)$ and $\abs{m_0(\xi)}$ for the optimal choice of the
root. Note that the Cohen criterion is clearly satisfied.

\begin{figure}
\begin{center}
\scalebox{0.6}{\includegraphics{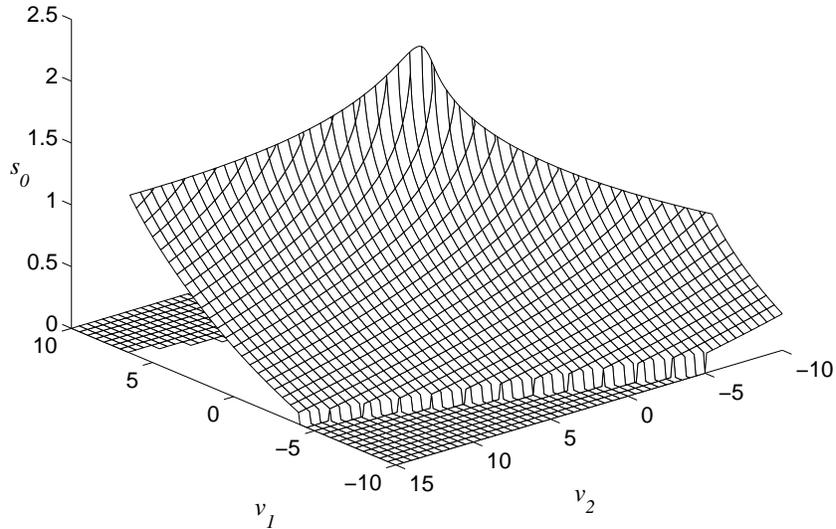}}
\end{center}
\caption{The Sobolev regularity exponent $s_0$ as a function of the
parameters $\vbf=(v_1,v_2)$ for the case $2N=10$. (When $\vbf$ is outside
the admissible region $\min r_\vbf\geq0$ we plot $s_0$ as zero.)}\mylabel{f:s0102}
\end{figure}

\begin{figure}
\begin{center}
\scalebox{0.6}{\includegraphics{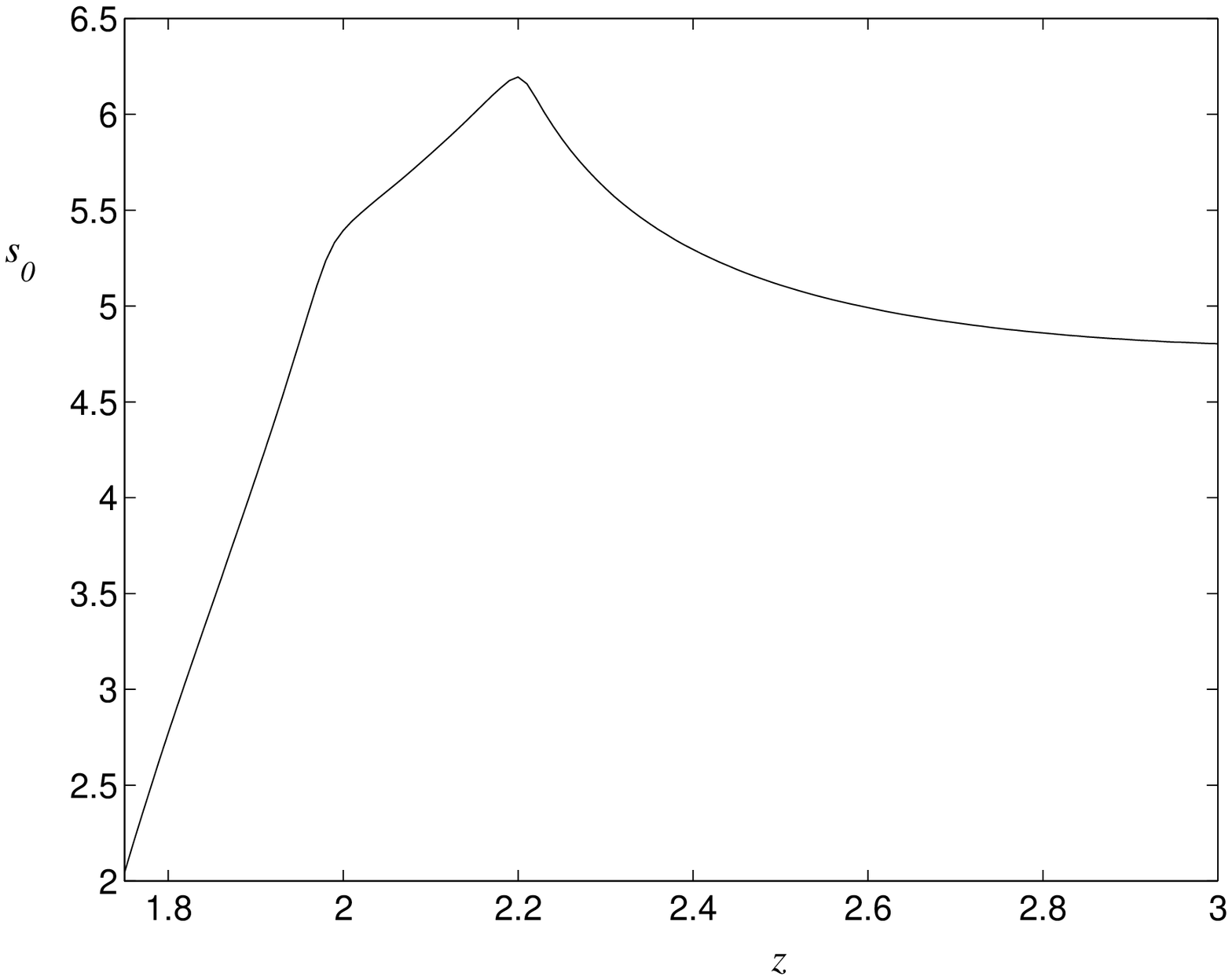}}
\end{center}
\caption{The Sobolev regularity exponent $s_0$ as a function of the
root $z$ that $m_0$ and $r$ have on the interval $\hpipi$. Here
$2N=32$, $\nz=1$.}\mylabel{f:s0z1}
\end{figure}

\begin{figure}
\begin{center}
\scalebox{0.6}{\includegraphics{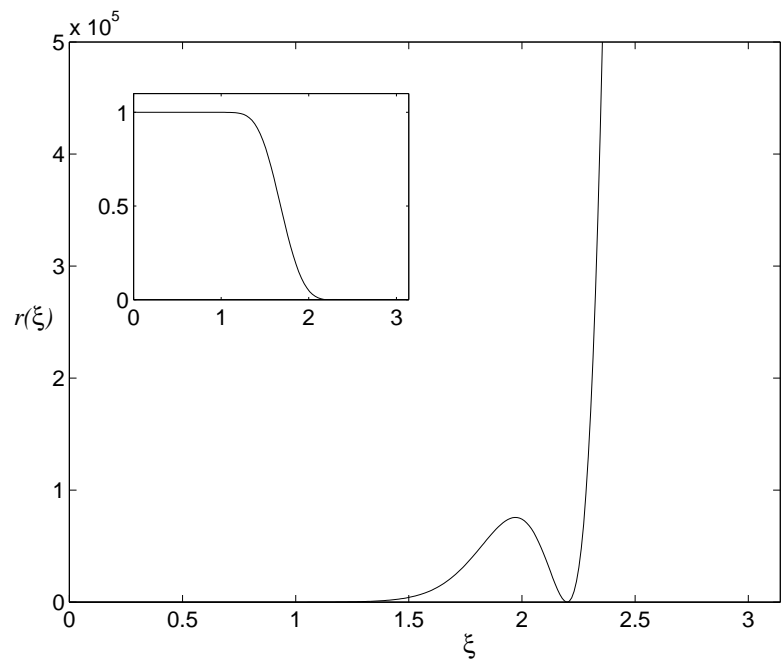}}
\end{center}
\caption{The function $r(\xi)$ for $2N=32$, $\nz=1$ corresponding to
optimal regularity.
The inset shows the graph of  $\abs{m_0(\xi)}$.}\mylabel{f:r322}
\end{figure}

Numerical optimization of the regularity exponent was done by using
the Matlab Optimization Toolbox.
The smoothest wavelet was always found on the boundary of the feasible
region $\{\vbf\in\Rbb^2\colon\min r_\vbf\geq0\}$, that is, the optimal
$r$ has a double root on $\hpipi$. 
For comparison with the Daubechies family, the most regular scaling
function with ten coefficients is shown in Figure~\ref{f:diffs}
together with its first derivative and graphs of the corresponding
Daubechies scaling function. Figure~\ref{f:m0prod} shows graphs of the
Fourier transforms of the scaling functions. 

\begin{figure}
\begin{center}
\begin{tabular}{rp{.02cm}r}
\scalebox{.375}{\includegraphics{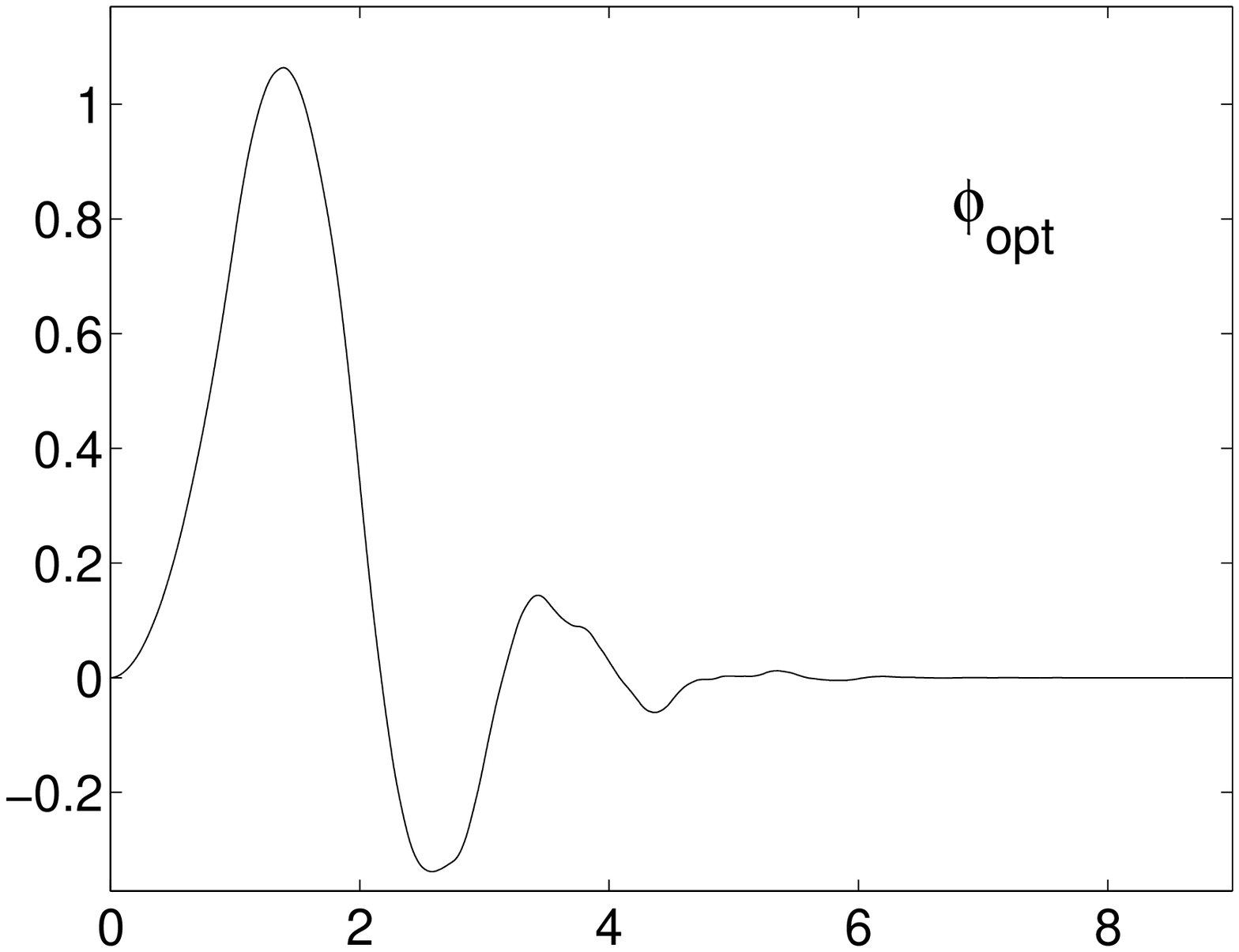}} &&
\scalebox{.375}{\includegraphics{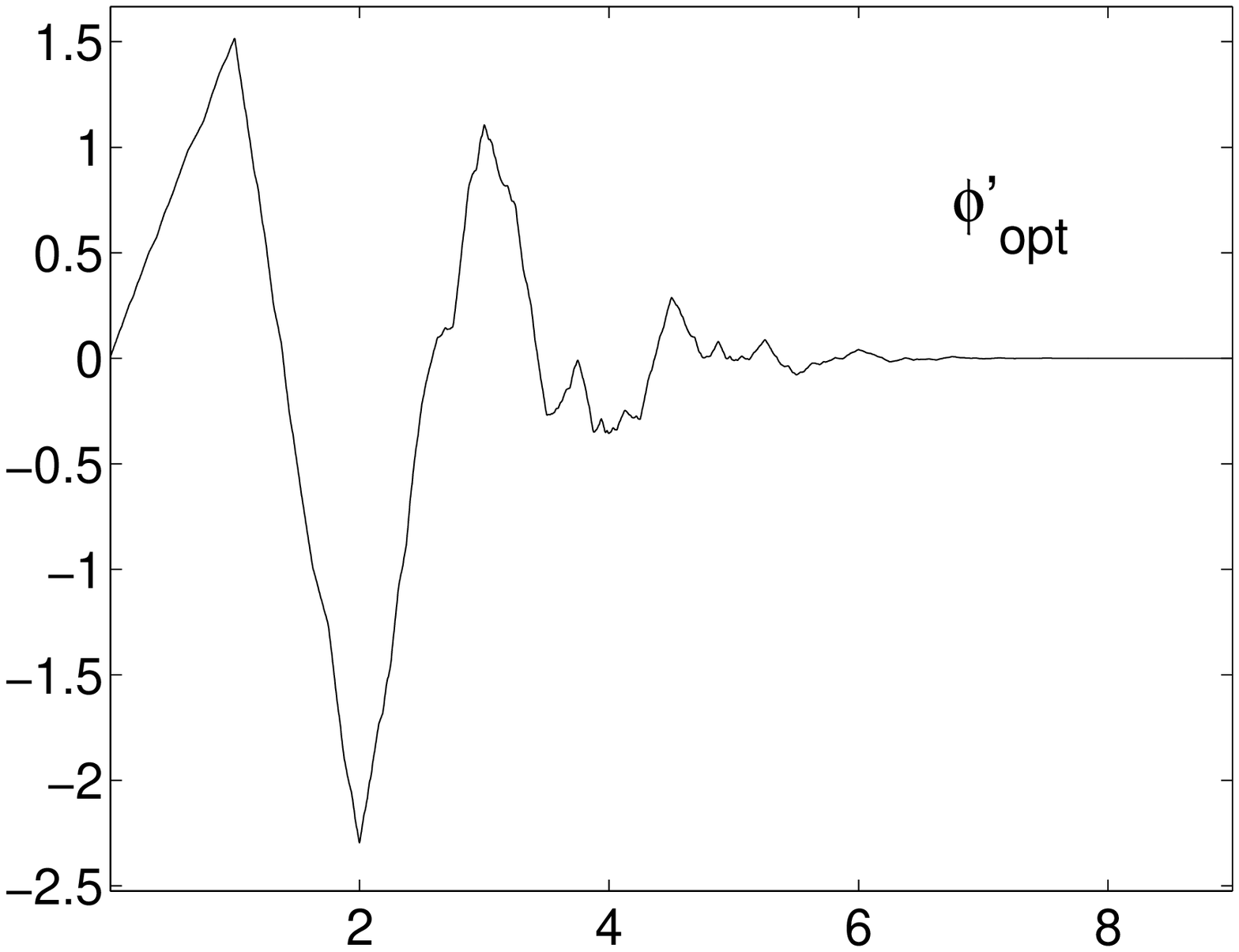}} \\
\scalebox{.375}{\includegraphics{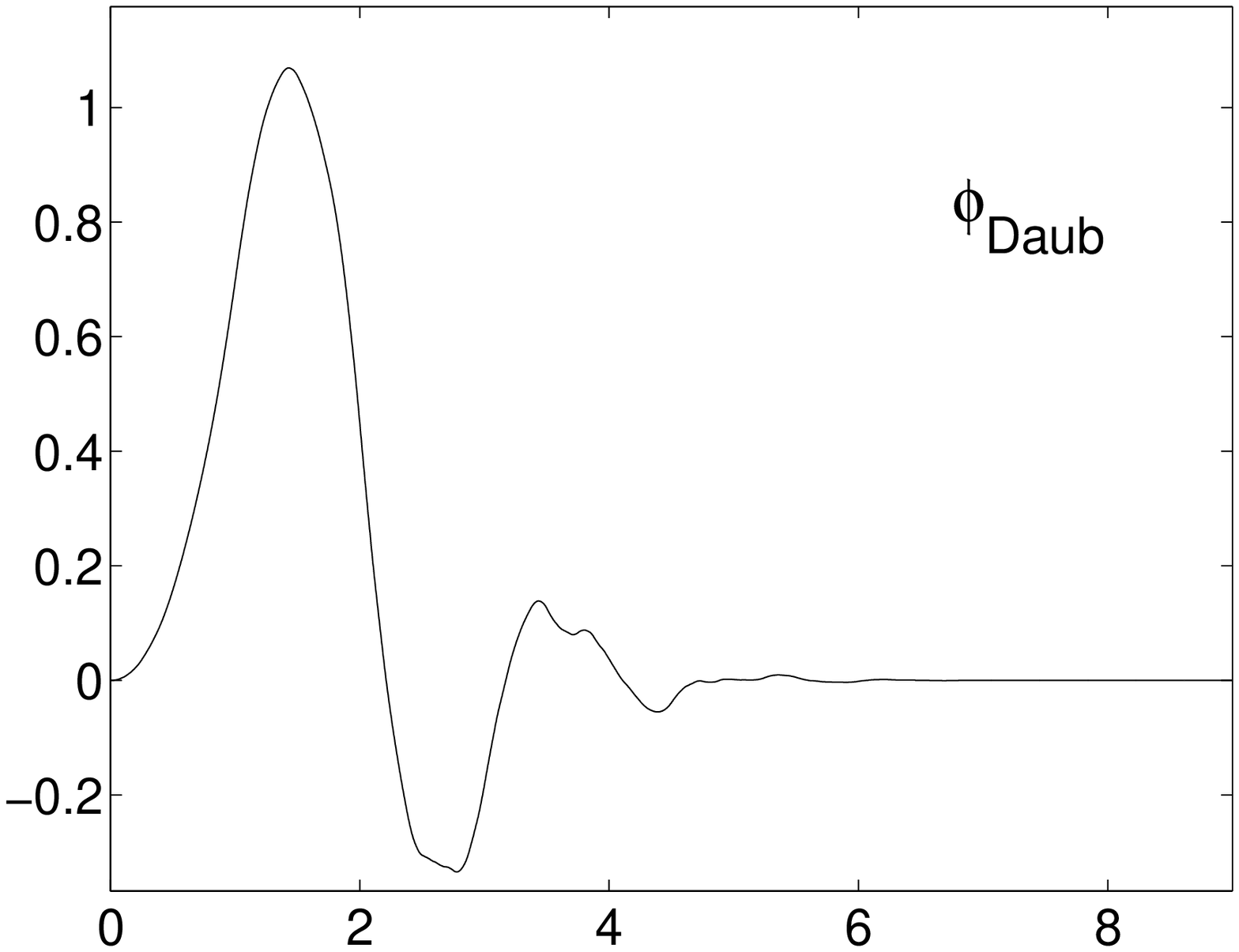}} &&
\scalebox{.375}{\includegraphics{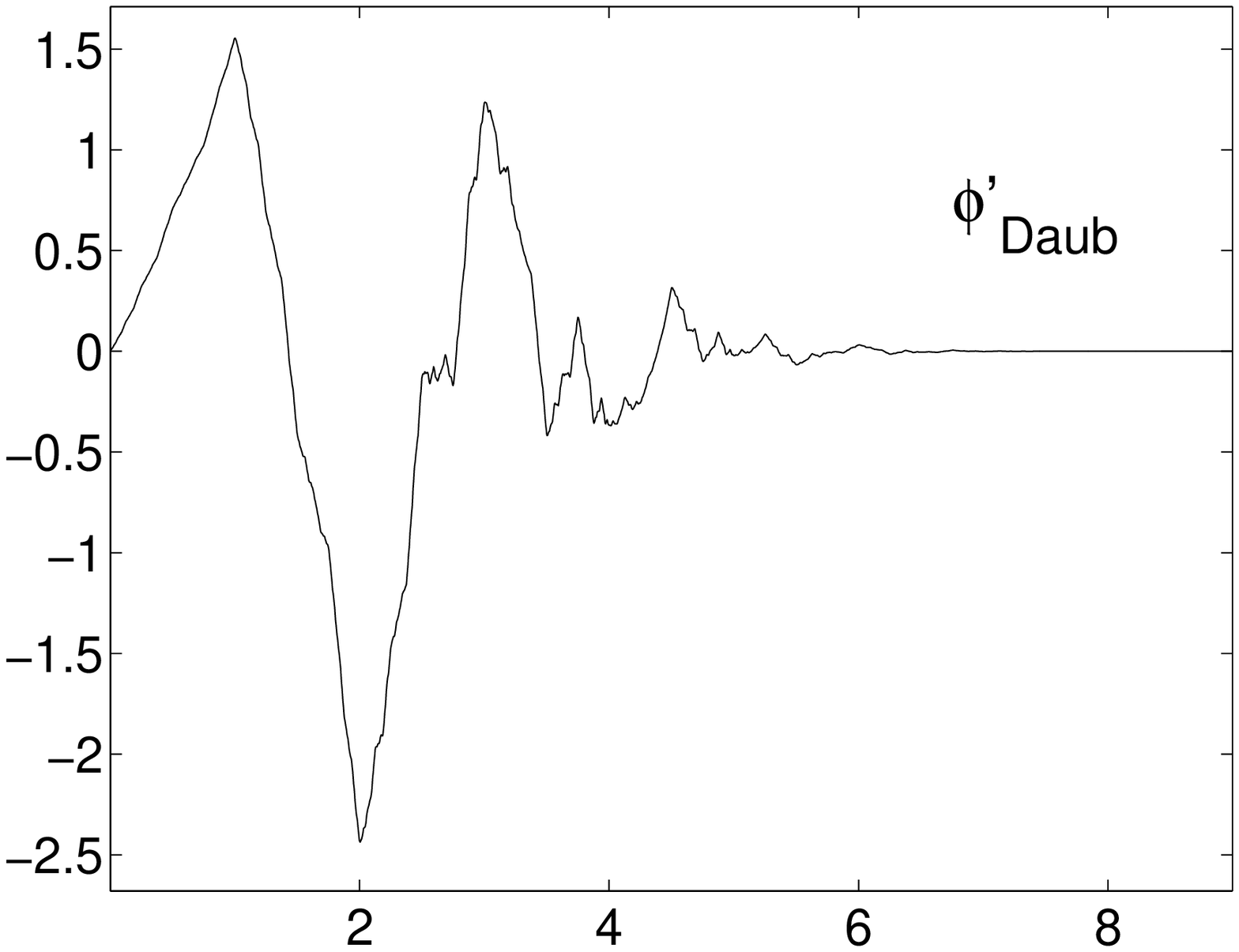}} \\
\end{tabular}
\end{center}
\caption{On the top the most regular scaling function and its first
derivative, on the bottom the Daubechies
scaling function and its first derivative. Here $2N=10$.}\mylabel{f:diffs}
\end{figure}

\begin{figure}
\begin{center}
\scalebox{0.6}{\includegraphics{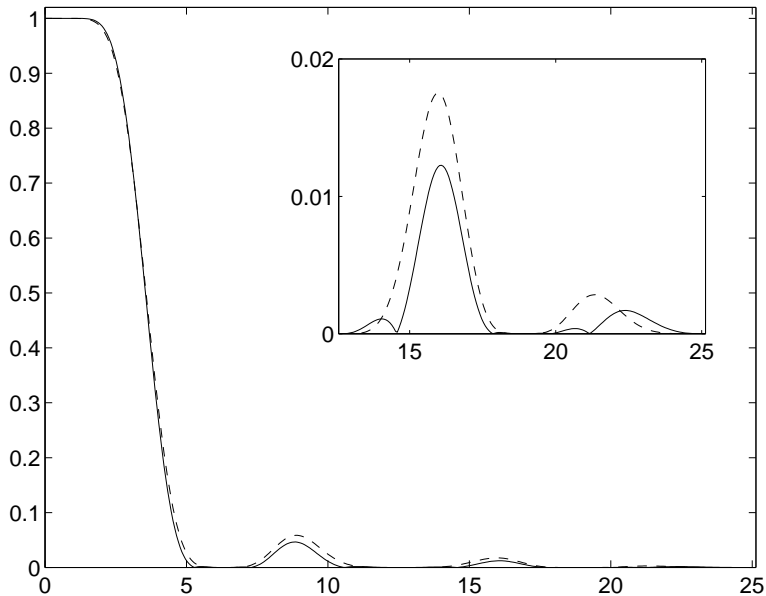}}
\end{center}
\caption{Graphs of $\widehat{\phi}(\xi)$ given by~\eqref{e:m0prod}
    for the most
    regular (solid curve) and Daubechies wavelets (dashed) with ten
    coefficients. The blow-up of the interval $[4\pi,8\pi]$ in the inset
    shows how the root of $m_0$ on $\hpipi$ improves decay 
    of~$\widehat\phi$, hence regularity of~$\phi$.
}\mylabel{f:m0prod}
\end{figure}

Table~\ref{t:z1} lists the locations of the root of $m_0$ on $\hpipi$
for the most regular wavelets.

\begin{table}
\begin{center}
\begin{tabular}{cp{1em}r@{.}lp{2cm}cp{1em}r@{.}l}
    $2N$  &&  \multicolumn{2}{c}{$z$}
                       &&    $2N$  &&  \multicolumn{2}{c}{$z$} \\[.5ex]
          \cline{1-4}\cline{6-9} \\[-2ex]
      8   &&   2&8762  &&     26   &&   2&2346  \\
     10   &&   2&6450  &&     28   &&   2&2213  \\
     12   &&   2&5099  &&     30   &&   2&2099  \\
     14   &&   2&4200  &&     32   &&   2&1995  \\
     16   &&   2&3630  &&     34   &&   2&1894  \\
     18   &&   2&3238  &&     36   &&   2&1799  \\
     20   &&   2&2939  &&     38   &&   2&1712  \\
     22   &&   2&2701  &&     40   &&   2&1633  \\
     24   &&   2&2506  \\
\end{tabular}
\end{center}
\caption{Location of the root of $r$ for the most regular wavelets
for $\nz=1$.}\mylabel{t:z1}
\end{table}

The set $\{\vbf\in\Rbb^2\colon\min r_\vbf\geq0\}$ contains a
parametrization of all wavelets with $2N$ coefficients and $N-2$
vanishing moments. Together with the results of the previous section
this shows the regularity exponents in the $\nz=1$ column of
Table~\ref{t:s0res} are the best possible among wavelets with at least
$N-2$ vanishing moments and $2N\geq8$ filter coefficients.

\subsection{More degrees of freedom}\mylabel{s:dofn}

The approach of the previous section did not work with higher degrees of
freedom, i.e., with $\vbf\in\Rbb^n$, $n\geq3$, due to convergence
problems with Matlab's constrained optimizer (since $r$ is a
polynomial of high degree, problems with the constraint $\min
r_\vbf\geq0$ are not unexpected). Instead, the optimization
was done by considering the Sobolev regularity exponent $s_0$ as a
function of the zeros $\zbf=(z_1,\dots,z_\nz)$ of $r$ on $\hpipi$. The
benefits are that now we are studying an unconstrained optimization
problem and the number of degrees of freedom is halved.

For a fixed value of $\zbf$ we can solve for the parameters $\vbf$
discussed above by requiring $r_\vbf(\zbf)=r'_\vbf(\zbf)=0$. This is a
linear system, though it may be badly conditioned when $r$ is of high
degree. (This problem was overcome by using Mathematica and 30 digit
precision in the computations.)

The optimal locations for the roots of $m_0$ are given in
Tables~\ref{t:z2},~\ref{t:z3}, and~\ref{t:z4}. Figure~\ref{f:s0z2}
shows a typical example of the regularity exponent as a function of
the roots and Figure~\ref{f:r328} shows a plot of $r(\xi)$ for the
optimal choice of the roots.  Again the Cohen criterion is clearly
satisfied.

As discussed in Section~\ref{s:res} it is still open what happens if
$r$ is not required to have as many roots on $\hpipi$ as above. Recall
that $r$ is essentially a polynomial (of cosines) of high
degree. Figure~\ref{f:r328} shows high peaks between the roots of
$r$. Is it possible to make these peaks lower by not requiring $r$ to
vanish between them? However, even if this were the case, it is not
clear how the regularity exponent would change.

\begin{table}
\begin{center}
\begin{tabular}{cp{1cm}r@{.}lp{1em}r@{.}l}
    $2N$  &&  \multicolumn{2}{c}{$z_1$} &&  \multicolumn{2}{c}{$z_2$} \\[.5ex] \hline\\[-2ex]
     16  &&  2&3525  &&  2&8336  \\
     18  &&  2&3150  &&  2&7491  \\
     20  &&  2&2790  &&  2&7110  \\
     22  &&  2&2496  &&  2&6867  \\
     24  &&  2&2260  &&  2&6691  \\
     26  &&  2&2072  &&  2&6566  \\
     28  &&  2&1931  &&  2&6448  \\
     30  &&  2&1814  &&  2&6337  \\
     32  &&  2&1718  &&  2&6255  \\
     34  &&  2&1630  &&  2&6217  \\
     36  &&  2&1557  &&  2&6199  \\
     38  &&  2&1501  &&  2&6192  \\
     40  &&  2&1450  &&  2&6186  
\end{tabular}
\end{center}
\caption{Location of the roots of $r$ for the most regular wavelets
for $\nz=2$.}\mylabel{t:z2}
\end{table}

\begin{table}
\begin{center}
\begin{tabular}{cp{1cm}r@{.}lp{1em}r@{.}lp{1em}r@{.}l}
    $2N$  &&  \multicolumn{2}{c}{$z_1$} &&  \multicolumn{2}{c}{$z_2$}
          &&  \multicolumn{2}{c}{$z_3$} \\[.5ex] \hline\\[-2ex]
     22  &&  2&2474  &&  2&6807  &&  3&0380  \\
     24  &&  2&2199  &&  2&6571  &&  2&9637  \\
     26  &&  2&2005  &&  2&6413  &&  2&9197  \\
     28  &&  2&1871  &&  2&6303  &&  2&8892  \\
     30  &&  2&1761  &&  2&6229  &&  2&8710  \\
     32  &&  2&1668  &&  2&6196  &&  2&8598  \\
     34  &&  2&1582  &&  2&6184  &&  2&8440  \\
     36  &&  2&1497  &&  2&6172  &&  2&8071  \\
     38  &&  2&1425  &&  2&6143  &&  2&7488  \\
     40  &&  2&1369  &&  2&6174  &&  2&7109  
\end{tabular}
\end{center}
\caption{Location of the roots of $r$ for the most regular wavelets
for $\nz=3$.}\mylabel{t:z3}
\end{table}

\begin{table}
\begin{center}
\begin{tabular}{cp{1cm}r@{.}lp{1em}r@{.}lp{1em}r@{.}lp{1em}r@{.}l}
    $2N$  &&  \multicolumn{2}{c}{$z_1$} &&  \multicolumn{2}{c}{$z_2$}
          &&  \multicolumn{2}{c}{$z_3$} &&  \multicolumn{2}{c}{$z_4$} \\[.5ex] \hline\\[-2ex]
     28  &&  2&2307  &&  2&4928  &&  2&9401  &&  3&0099  \\
     30  &&  2&1767  &&  2&6544  &&  2&6958  &&  3&0326  \\
     32  &&  2&1657  &&  2&6250  &&  2&8334  &&  2&9133  \\
     34  &&  2&1581  &&  2&5655  &&  2&6837  &&  2&9097  \\
     36  &&  2&1480  &&  2&5568  &&  2&6178  &&  2&9207  \\
     38  &&  2&1362  &&  2&4881  &&  2&6498  &&  2&9159  \\
     40  &&  2&1316  &&  2&6311  &&  2&7035  &&  2&8414  
\end{tabular}
\end{center}
\caption{Location of the roots of $r$ for the most regular wavelets
for $\nz=4$.}\mylabel{t:z4}
\end{table}

\begin{figure}
\begin{center}
\scalebox{0.6}{\includegraphics{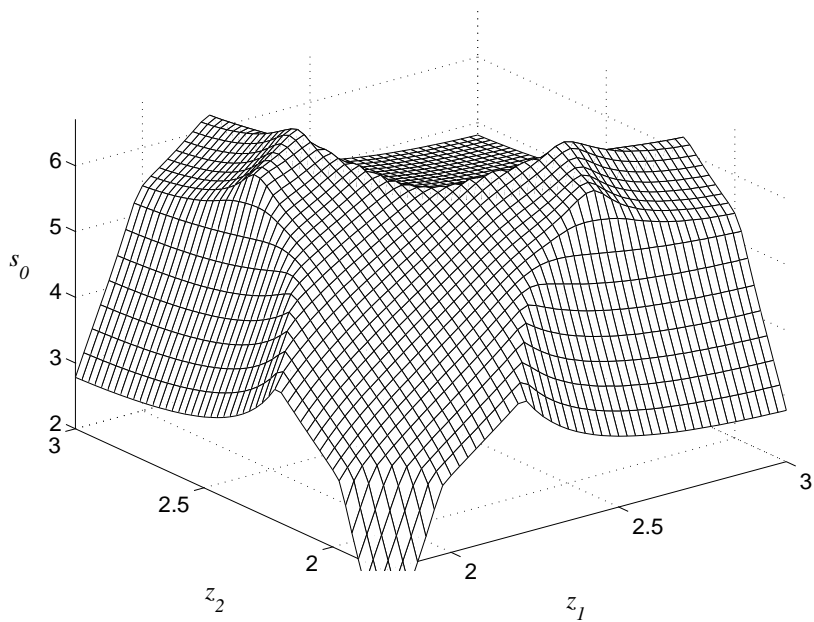}}
\end{center}
\caption{The Sobolev regularity exponent $s_0$ as a function of the
roots $z_1$, $z_2$ that $m_0$ and $r$ have on the interval
$\hpipi$. Here $2N=32$.}
\mylabel{f:s0z2}
\end{figure}

\begin{figure}
\begin{center}
\scalebox{0.6}{\includegraphics{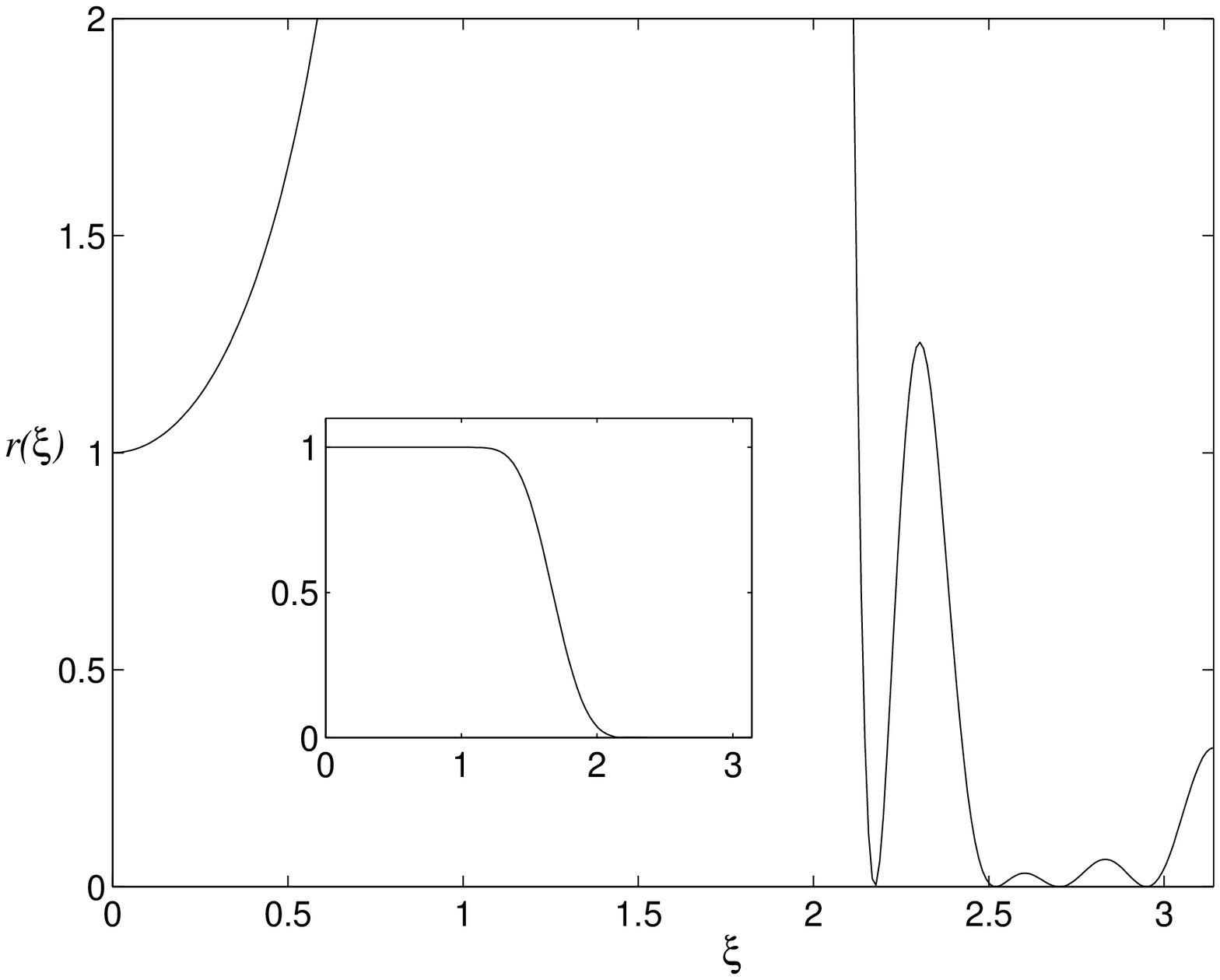}}
\end{center}
\caption{The function $r(\xi)$ for $2N=32$, $\nz=4$ corresponding to
optimal regularity.
The inset shows the graph of $\abs{m_0(\xi)}$.}\mylabel{f:r328}
\end{figure}

\section{Other parametrizations}\mylabel{s:other-param}

The parametrization of the wavelets could also be based on the
original result of Daubechies \cite[Proposition~6.1.2]{Daub92}, or on
a similar result by Wells~\cite{Well93}. Both of these were tried, but
they suffer from the same difficulty as our parametrization in terms
of the vector $\vbf$: there is a non-linear constraint that makes
numerical optimization for complicated filters impossible with the
software that was used. (These parametrizations could be used as the
starting point for a parametrization in terms of roots of $m_0$,
though.)

Direct parametrizations of the coefficients $c_k$ have been introduced
in~\cite{Lina93,Poll90,Schn93,Zou93}, but in contrast to our
approach the moment condition~\eqref{e:m0def} is then a system of
non-linear equations of the parameters (see in particular
\cite[Theorem~2]{Zou93}). For wavelets with $2N$ filter coefficients and
$M$ vanishing moments there are $N$ parameters and $M-1$ non-linear
equations. The equations can be used to solve for $M-1$ of the
parameters when the rest are kept fixed. Numerical
experiments show, however, that typically there are multiple
solutions that correspond to wavelets with \emph{different regularity},
making this approach unsuitable for optimizing regularity. Moreover, these
parametrizations provide an expression directly for $m_0(\xi)$, but for the
computation of the regularity exponent we need $r(\xi)$.

Only the parametrization in terms of the roots of $m_0$ has provided
the computational simplicity necessary for optimizing regularity when
the filter is complicated: recall that there are only linear equations
to be solved.

\section{Filter coefficients}\mylabel{s:coef}

The optimization provides us with an expression for $\abs{m_0(\xi)}$.
The (low-pass) filter coefficients $c_k$ in~\eqref{e:dileqn} are then
easy to compute by using spectral factorization
\cite[p.~172]{Daub92}. The coefficients for selected most
regular wavelets are listed in Tables
\ref{t:fc101}--\ref{t:fc303}.  To guarantee that $m_0$ satisfies the
orthonormality condition~\eqref{e:ortho} to high precision, the
factorization was done with Mathematica using 30 digit precision and
checked by repeating the computations with 60 digit precision: the
results are accurate to the number of decimals given in the tables.

\begingroup

\newenvironment{wfc}[3]
  {\begingroup
    \renewcommand{\baselinestretch}{.97}
    \newcommand{\capt}{$2N=#1$, $\nz=#2$}
    \newcommand{\lbl}{t:fc#1#2}
    \begin{table}[p]
    \begin{center}\tt
    \begin{tabular}{r@{.}l}
    \multicolumn{2}{c}{$c_k$} \\[.5ex]\hline
	\multicolumn{2}{c}{}\\[-2ex]}
  {\end{tabular}
    \end{center}
    \caption{\capt}\mylabel{\lbl}
    \end{table}
 \endgroup}

  \begin{wfc}{10}{1}{3}
     1&807084186243315e-1\\
     6&272955371644549e-1\\
     7&021176047824324e-1\\
     1&120128480002895e-1\\
    -2&446469866533168e-1\\
    -2&970511286791226e-2\\
     8&518480911807988e-2\\
    -7&179751673520435e-3\\
    -1&625706468497944e-2\\
     4&683260563235789e-3
  \end{wfc}

  \begin{wfc}{20}{2}{6}
     4&409469394058257e-2\\
     2&599093138855805e-1\\
     6&057694016921036e-1\\
     6&371421499968423e-1\\
     1&364367462665367e-1\\
    -2&883826797756209e-1\\
    -1&270531514862528e-1\\
     1&555314597190144e-1\\
     6&802158361295070e-2\\
    -8&862450100620332e-2\\
    -2&399047021783915e-2\\
     4&548910378211001e-2\\
     2&025062066984338e-3\\
    -1&790050908684461e-2\\
     3&315675665790194e-3\\
     4&335970768455452e-3\\
    -1&852310988320677e-3\\
    -3&359209223085526e-4\\
     3&395506340120817e-4\\
    -5&760617447778257e-5
  \end{wfc}

  \begin{wfc}{30}{3}{9}
     9&286962261105670e-3\\
     8&065009178408755e-2\\
     2&975638977797063e-1\\
     5&849583107003815e-1\\
     5&940950288790657e-1\\
     1&464131039850201e-1\\
    -2&899897992799635e-1\\
    -1&888671876256031e-1\\
     1&529318071875427e-1\\
     1&393659002859078e-1\\
    -9&487951734281382e-2\\
    -8&688845136891955e-2\\
     6&416864747737246e-2\\
     4&586586496411027e-2\\
    -4&242829324588881e-2\\
    -1&870129886474555e-2\\
     2&492219541071851e-2\\
     4&323699937678751e-3\\
    -1&195215551785778e-2\\
     7&996165686836154e-4\\
     4&259303332169016e-3\\
    -1&263090636020650e-3\\
    -9&623632092949613e-4\\
     5&655466161125865e-4\\
     7&527472879191102e-5\\
    -1&260995097154413e-4\\
     2&020452355886283e-5\\
     1&026632536894423e-5\\
    -4&411797664714925e-6\\
     5&080242007105260e-7
  \end{wfc}

\endgroup

\section{Conclusions}

We have constructed wavelets that are more regular than any of the wavelets
that have appeared in the literature with the same number of filter
coefficients, thus improving the ratio of regularity to complexity.
Our methods show that a numerical optimization approach provides the
flexibility necessary to construct highly regular wavelets while
keeping the number of vanishing moments of the wavelet as high as
possible.

We have also shown our results are optimal for wavelets with
$2N$~filter coefficients and at least $N-2$~vanishing moments.

\bibliographystyle{amsplain}
\bibliography{mrabbrev,analysis}

\providecommand{\bysame}{\leavevmode\hbox to3em{\hrulefill}\thinspace}
\begin{thebibliography}{10}

\bibitem{Cohe90}
A.~Cohen, \emph{Ondelettes, analyses multir\'esolutions et filtres miroirs en
  quadrature}, Ann. Inst. H. Poincar\'e Anal. Non Lin\'eaire \textbf{7} (1990),
  no.~5, 439--459.

\bibitem{Cohe92}
A.~Cohen and J.-P. Conze, \emph{R\'egularit\'e des bases d'ondelettes et
  mesures ergodiques}, Rev. Mat. Iberoamericana \textbf{8} (1992), no.~3,
  351--365.

\bibitem{Daub88}
I.~Daubechies, \emph{Orthonormal bases of compactly supported wavelets}, Comm.
  Pure Appl. Math. \textbf{41} (1988), 909--996.

\bibitem{Daub92}
\bysame, \emph{Ten Lectures on Wavelets}, Society for Industrial and Appl.
  Math., 1992.

\bibitem{Daub93}
\bysame, \emph{Orthonormal bases of compactly supported wavelets {I}{I}.
  {V}ariations on a theme}, SIAM J. Math. Anal. \textbf{24} (1993), no.~2,
  499--519.

\bibitem{Daub92a}
I.~Daubechies and J.~C. Lagarias, \emph{Two-scale difference equations. {I}{I}.
  {L}ocal regularity, infinite products of matrices and fractals}, SIAM J.
  Math. Anal. \textbf{23} (1992), no.~4, 1031--1079.

\bibitem{Eiro92}
T.~Eirola, \emph{Sobolev characterization of solutions of dilation equations},
  SIAM J. Math. Anal. \textbf{23} (1992), no.~4, 1015--1030.

\bibitem{Lema98}
P.~G. Lemari\'e-Rieusset and E.~Zahrouni, \emph{More regular wavelets}, Appl.
  Comput. Harmon. Anal. \textbf{5} (1998), 92--105.

\bibitem{Lina93}
J.-M. Lina and M.~Mayrand, \emph{Parametrizations for {D}aubechies wavelets},
  Phys. Rev. E (3) \textbf{48} (1993), no.~6, R4160--R4163.

\bibitem{Ojan91}
H.~Ojanen, \emph{Remarks on the {S}obolev regularity of wavelets and
  interpolation schemes}, Research Reports A305, Helsinki Univ. of Technology,
  Inst. of Math., 1991.

\bibitem{Poll90}
D.~Pollen, \emph{${\rm {S}{U}}\sb {I}(2,{F}[z,1/z])$ for ${F}$ a subfield of
  ${\bf {C}}$}, J. Amer. Math. Soc. \textbf{3} (1990), no.~3, 611--624.

\bibitem{Schn93}
J.~Schneid and S.~Pittner, \emph{On the parametrization of the coefficients of
  dilation equations for compactly supported wavelets}, Computing \textbf{51}
  (1993), no.~2, 165--173.

\bibitem{Vill92}
L.~F. Villemoes, \emph{Energy moments in time and frequency for two-scale
  difference equation solutions and wavelets}, SIAM J. Math. Anal. \textbf{23}
  (1992), no.~6, 1519--1543.

\bibitem{Volk92}
H.~Volkmer, \emph{On the regularity of wavelets}, IEEE Trans. Inform. Theory
  \textbf{38} (1992), no.~2, part 2, 872--876.

\bibitem{Volk95}
\bysame, \emph{Asymptotic regularity of compactly supported wavelets}, SIAM J.
  Math. Anal. \textbf{26} (1995), no.~4, 1075--1087.

\bibitem{Well93}
Wells, R.~O., Jr., \emph{Parametrizing smooth compactly supported wavelets},
  Trans. Amer. Math. Soc. \textbf{338} (1993), no.~2, 919--931.

\bibitem{Zou93}
H.~Zou and A.~H. Tewfik, \emph{Parametrization of compactly supported
  orthonormal wavelets}, IEEE Trans. Signal Processing \textbf{41} (1993),
  no.~3, 1428--1431.

\end{thebibliography}

\end{document}